\documentclass[11pt]{article}

\usepackage{enumerate}
\usepackage{amsmath,xspace,amssymb,mathrsfs}

\input xy
\xyoption{all}
\xyoption{2cell}
\UseAllTwocells
\CompileMatrices

\title{Limits of small functors \thanks{Both authors gratefully acknowledge the
support of the Australian Research Council.} }
\author{Brian J. Day \\
Mathematics Department\\
Macquarie University\\
\and
Stephen Lack \\
School of Computing and Mathematics \\
University of Western Sydney \\
{\tt s.lack@uws.edu.au}
}
\date{}

\renewcommand{\phi}{\varphi}

\newcommand{\A}{{\ensuremath{\mathscr{A}}}\xspace}
\newcommand{\B}{{\ensuremath{\mathscr{B}}}\xspace}
\newcommand{\BB}{{\ensuremath{\overline{\mathscr{B}}}}\xspace}
\newcommand{\C}{{\ensuremath{\mathscr{C}}}\xspace}
\newcommand{\D}{{\ensuremath{\mathscr{D}}}\xspace}
\newcommand{\DD}{{\ensuremath{\mathscr{D'}}}\xspace}
\newcommand{\E}{{\ensuremath{\mathscr{E}}}\xspace}
\newcommand{\G}{{\ensuremath{\mathscr{G}}}\xspace}
\newcommand{\I}{{\ensuremath{\mathscr{I}}}\xspace}
\newcommand{\K}{{\ensuremath{\mathscr{K}}}\xspace}
\newcommand{\LL}{{\ensuremath{\mathscr{L}}}\xspace}
\newcommand{\M}{{\ensuremath{\mathscr{M}}}\xspace}
\newcommand{\N}{{\ensuremath{\mathbb{N}}}\xspace}
\renewcommand{\P}{{\ensuremath{\mathscr{P}}}\xspace}
\newcommand{\V}{{\ensuremath{\mathscr{V}}}\xspace}

\newcommand{\PC}{\ensuremath{\mathscr{P\!C}}\xspace}
\newcommand{\PK}{\ensuremath{\mathscr{P\!K}}\xspace}
\newcommand{\PL}{\ensuremath{\mathscr{P\!L}}\xspace}

\newcommand{\Set}{\textnormal{\bf Set}\xspace}
\newcommand{\SSet}{\textnormal{\bf SSet}\xspace}
\newcommand{\Ab}{\textnormal{\bf Ab}\xspace}
\newcommand{\RMod}{\textnormal{\bf $R$-Mod}\xspace}
\newcommand{\Gpd}{\textnormal{\bf Gpd}\xspace}
\newcommand{\Cat}{\textnormal{\bf Cat}\xspace}
\newcommand{\Ban}{\textnormal{\bf Ban}\xspace}
\newcommand{\VCat}{\textnormal{\bf \V-Cat}\xspace}

\newcommand{\op}{\ensuremath{^{\textnormal{op}}}}
\newcommand{\ev}{\textnormal{ev}}
\newcommand{\Lan}{\textnormal{Lan}}
\newcommand{\colim}{\textnormal{colim}}
\newcommand{\cone}{\textnormal{cone}}
\newcommand{\ot}{\otimes}

\newcommand{\ct}{\pitchfork}

\newdir{ >>}{{}*!/-10pt/@{>>}}
\newdir{ >}{{}*!/-8pt/@{>}}

\renewcommand{\t}{\times}         
\def\endproof{{\parfillskip=0pt\hfill$\Box$\vskip 10pt}}

\newcommand{\two}{\ensuremath{{\hbox{\textrm 2}\kern-.25em
        \hbox{\vrule height1.5ex width 0.4pt depth -.2ex}}\kern.2em}\xspace}

\newtheorem{theorem}{Theorem}[section]    
\newtheorem{corollary}[theorem]{Corollary}   
\newtheorem{proposition}[theorem]{Proposition}   
\newtheorem{lemma}[theorem]{Lemma}   
   
\newtheorem{preremark}[theorem]{Remark}   
\newtheorem{prexample}[theorem]{Example}   

\newenvironment{remark}{\begin{preremark}\rm}{\end{preremark}}
\newenvironment{example}{\begin{prexample}\rm}{\end{prexample}}

\newcommand{\proof}{\noindent{\sc Proof:}\xspace}

\begin{document}

\label{firstpage}
\maketitle

\begin{abstract}
For a small category \K enriched over a suitable monoidal category \V,
the free completion of \K under colimits is the presheaf category 
$[\K\op,\V]$. If \K is large, its free completion under colimits is
the \V-category \PK of small presheaves on \K, where a presheaf is 
small if it is a left Kan extension of some presheaf with small
domain. We study the existence of limits and of monoidal closed 
structures on \PK.
\end{abstract}

A fundamental construction in category theory is the category of 
presheaves $[\K\op,\Set]$ on a small category \K. Among many other
important properties, it is the free completion of \K under colimits.
If the category \K is large, then the full presheaf category $[\K\op,\Set]$
is not the free completion of \K under colimits; indeed it is not even
a legitimate category, insofar as its hom-sets are not in general small.

In some contexts it is more appropriate to consider not {\em all} the 
presheaves on \K, but only the {\em small} ones: a presheaf
$F:\K\op\to\Set$ is said to be small if it is the left Kan extension
of some presheaf whose domain
is small. This is equivalent to $F$ being the left Kan extension of its
restriction to some small full subcategory of its domain, or equally 
to its being a 
small colimit of representables. The natural transformations between 
two small presheaves on \K do form a small set, and so the totality of
small presheaves on \K forms a genuine category \PK with small hom-sets. 
Furthermore, \PK is in fact the free completion of \K under colimits.
Of course if \K is small, then every presheaf on \K is small, and so 
\PK is just $[\K\op,\Set]$, but in general this is not the case. 

Although \PK is the free completion of \K under colimits, it does not
have all the good properties of $[\K\op,\Set]$ for small \K. For example
it is not necessarily complete or cartesian closed. In this paper we
study, among other things, when \PK does have such good properties.

In fact we work not just with ordinary categories, but with categories
enriched over a suitable monoidal category \V. Once again, if \K is 
small then $[\K\op,\V]$ is the free completion of \K under colimits,
but for large \K this is no longer the case; the illegitimacy of 
$[\K\op,\V]$ in that case is more drastic: it is not even a \V-category.
The free completion of \K under colimits is the \V-category \PK of 
small presheaves on \K, where once again a presheaf is small if it is 
the left Kan extension of some presheaf with small domain; and once 
again the two reformulations of this notion can be made.

The case $\V=\Set$ is closely related to work by various authors. Freyd
\cite{Freyd-lucid} introduced two smallness notions for presheaves on large
categories. He called a functor $F:\K\op\to\Set$ {\em petty} if there is a
small family $(C_\lambda\in\K)_{\lambda\in\Lambda}$ with an epimorphism 
$$\sum_{\lambda}\K(-,C_{\lambda})\to F;$$
and {\em lucid} if it is petty and for any representable $\K(-,A)$ and
any pair of maps $u,v:\K(-,A)\to F$, their equalizer is petty. 
%(Our
%notion of {\em weak smallness}, studied in Section~\ref{sect:P_w} below,
%is a generalization to the enriched setting of Freyd's petty functors.)
Freyd studied when the category of petty presheaves on \K is complete, and when
the category of lucid presheaves on \K is complete, obtaining results similar
to our Theorem~\ref{thm:main} below. Rosick\'y \cite{Rosicky-ccc} showed
that if \K is complete, then a presheaf $F$ on \K is lucid if and only
if it is small; one can then deduce our Corollary~\ref{cor:complete} from
the results of Freyd. Rosick\'y also characterized, in the case $\V=\Set$,
when \PK is cartesian closed; see Example~\ref{ex:Rosicky} below.
In a slightly different direction, the existence of limits in free
completions under some class of colimits was studied in \cite{Karazeris-completeness}.

In the enriched case, the fact, mentioned above, 
that \PK is the free completion of \K under colimits, is due to 
Lindner~\cite{Lindner}. The existence of limits or monoidal closed
structures on \PK seems not to have been considered in the enriched
setting.

Some of our results have been used in abstract homotopy theory; for
example Corollary~\ref{cor:complete} was used in \cite{Chorny}. The idea
is that one wants to have a complete and cocomplete category of diagrams of
some particular type, where the indexing category is large. In this
context one is particularly interested in the case $\V=\SSet$, the
category of simplicial sets.

In Section~\ref{sect:review} we review the required background from 
enriched category theory, and in Section~\ref{sect:small} the notion of
small functor. Then in Section~\ref{sect:limits} we prove the fundamental
result that \PK is complete if and only if it has limits of representables;
thus in particular \PK is complete if \K is so. In Section~\ref{sect:sound}
we refine the results of the previous section to deal not with arbitrary
(small) limits, but with limits of some particular type, such as finite
limits or finite products. In Section~\ref{sect:ordinary} we deduce from 
the earlier results various known results about the case $\V_0=\Set$
of ordinary categories, before extending them to the case where $\V_0$
is a presheaf category.
%Section~\ref{sect:P_w} concerns weakly small functors, generalizing
%Freyd's petty functors to the enriched setting. 
Section~\ref{sect:continuity}
concerns not the existence of limits in \PK but the preservation of 
limits by functors $\P F:\PK\to\P\LL$ given by left Kan extension along
$F\op:\K\op\to\LL\op$. In Section~\ref{sect:promonoidal} we study 
monoidal closed structures on \PK using the notion of promonoidal 
category. In Section~\ref{sect:M} we consider limits of small 
functors with codomain a locally presentable category \M, generalizing
the earlier case of $\M=\V$. Finally in Section~\ref{sect:Isbell} we 
briefly discuss Isbell conjugacy for large categories.

The second-named author is very grateful to Francis Borceux and his
colleagues at the Universit\'e Catholique de Louvain-la-Neuve for their
hospitality during a one-month visit in 1998, during which some of the early
work on this paper was completed; and to Ross Street and the
Mathematics Department at Macquarie University, for their hospitality
during a sabbatical visit in 2006, during which the paper was finally
completed.

%%%%%%%%%%%%%%%%%%%%%%%%%%%%%%%%%%%%%%%%%%%%%%%%%%%%%%%%%%%%%%%%%%%%%%%%

\section{Review of relevant enriched category theory}
\label{sect:review}

%%%%%%%%%%%%%%%%%%%%%%%%%%%%%%%%%%%%%%%%%%%%%%%%%%%%%%%%%%%%%%%%%%%%%%%%

We shall work over a symmetric monoidal closed category \V. The
tensor product is denoted $\ot$, the unit object $I$, and the internal
hom $[~,~]$. Where 
necessary the underlying ordinary category is denoted $\V_0$. 

We suppose that this underlying ordinary category is locally presentable
\cite{Gabriel-Ulmer, AR}:
thus for some regular cardinal $\alpha$ and some small category \C with 
$\alpha$-small limits $\V_0$ is equivalent to the category of 
$\alpha$-continuous functors from \C to \Set. It follows that $\V_0$
is complete and cocomplete, and it turns out that \C is 
equivalent to the opposite of the full subcategory $(\V_0)_\alpha$
of $\V_0$ consisting of the $\alpha$-presentable objects: these are
the $X\in\V_0$ for which $\V_0(X,-):\V_0\to\Set$ preserves 
$\alpha$-filtered colimits. By \cite{vcat}, after possibly changing
$\alpha$, we may suppose that $(\V_0)_\alpha$ is closed in $\V_0$ 
under the monoidal structure, so that \V is {\em locally $\alpha$-presentable
as a closed category}, in the sense of \cite{Kelly-amiens}.

We shall work throughout the paper over such a locally presentable
closed category. This includes many important examples, such as the
categories 
\Set, \Ab, \RMod, \Cat, \Gpd, and \SSet, of sets, abelian groups, $R$-modules 
(over a commutative ring $R$), categories, groupoids, and simplicial sets,
as well as the two-element lattice \two. All these examples are locally
{\em finitely} presentable (that is, locally $\aleph_0$-presentable) but
there are further examples which require a higher cardinal than $\aleph_0$:
for example any Grothendieck topos, the category \Ban of Banach spaces
and linear contractions, Lawvere's category $[0,\infty]$ of 
extended non-negative real numbers, or the first-named author's
$*$-autonomous category $[-\infty,\infty]$ of extended real numbers. 
All categorical notions are 
understood to be enriched over \V, even if this is not explicitly stated.
(Thus category means \V-category, functor means \V-functor, and so on.)
We fix a regular cardinal $\alpha_0$ 
for which $\V_0$ is locally $\alpha_0$-presentable and $(\V_0)_{\alpha_0}$ 
is closed under the monoidal structure.
Henceforth ``$\alpha$ is a regular cardinal'' will mean ``$\alpha$ is 
a regular cardinal and $\alpha\ge\alpha_0$''.

For such a \V, it was shown in \cite{Kelly-amiens} that there is a good 
notion of locally $\alpha$-presentable 
\V-category, for any regular cardinal $\alpha\ge\alpha_0$. A locally 
$\alpha$-presentable \V-category \K is complete and cocomplete, and is 
equivalent to the \V-category of $\alpha$-continuous \V-functors from 
\C to \V for some small \V-category \C with $\alpha$-small limits. This 
\C can be identified with the opposite
of the category of $\alpha$-presentable objects in \K. 

%A weight $F:\K\op\to\V$ is said to be $\alpha$-flat if 
%$\alpha*-:\P(\K\op)\to\V$ is $\alpha$-continuous. Often we call
%a colimit $\alpha$-filtered if its weight is $\alpha$-flat.
%A \V-category \K is said to be $\alpha$-accessible if it is the free 
%completion under $\alpha$-filtered colimits of a small category
%\C. The locally $\alpha$-presentable categories can be characterized
%as the complete $\alpha$-accessible categories, or alternatively as the
%cocomplete $\alpha$-accessible categories (see \cite{BQR}). 
%A \V-category is accessible if it is $\alpha$-accessible for some $\alpha$,
%and locally presentable if it is locally $\alpha$-presentable for 
%some $\alpha$. (As usual we suppose that $\alpha\ge\alpha_0$.)

A {\em weight} is a presheaf $F:\C\op\to\V$, usually, although not always
with small domain. The {\em colimit} of a functor $S:\C\to\K$ is denoted
by $F*S$, while the {\em limit} of a functor $S:\C\op\to\K$ is denoted
by $\{F,S\}$. When $\C\op$ is the unit \V-category \I, we may identify
$F$ with an object of \V and $S$ with an object of \C; we sometimes 
write $F\cdot S$ for $F*S$ and call it a tensor, and we sometimes 
write $F\pitchfork S$ for $\{F,S\}$ and call it a cotensor.

%%%%%%%%%%%%%%%%%%%%%%%%%%%%%%%%%%%%%%%%%%%%%%%%%%%%%%%%%%%%%%%%%%%%%%%%

\section{Small functors}
\label{sect:small}

%%%%%%%%%%%%%%%%%%%%%%%%%%%%%%%%%%%%%%%%%%%%%%%%%%%%%%%%%%%%%%%%%%%%%%%%

A functor $F:\K\to\V$ is said to be {\em small} if it is the left Kan
extension of its restriction to some small full subcategory of \K. This
will clearly be the case if $F$ is a small colimit of representables, for
then we may take as the subcategory precisely those objects corresponding
to the representables in the colimit. On the other hand, if $F:\K\to\V$
is the left Kan extension of $FJ$ along the inclusion $J:\C\to\K$ of some
small full subcategory, then $F=(FJ)*\K(J,1)$, and so $F$ is a small 
colimit of representables. Thus the small functors are precisely the small
colimits of representables.

Of course if \K is itself small, then every functor from \K to \V is
small. If on the other hand \K is locally presentable, then a functor 
$F:\K\to\V$ is small if and only if it is accessible: that is, if and only 
if it preserves $\alpha$-filtered colimits for some regular cardinal 
$\alpha$. For if $F$ is accessible, then we may choose $\alpha$ so that 
\K is locally $\alpha$-presentable and $F$ preserves $\alpha$-filtered 
colimits; then $F$ is the left Kan extension of its restriction to the full
subcategory of \K consisting of the $\alpha$-presentable objects. 
Conversely, if $F$ is the left Kan extension of its restriction to a
small full subcategory \C of \K, then we may choose a regular cardinal 
$\alpha$ in such a way that \K is locally $\alpha$-presentable and every 
object in \C is $\alpha$-presentable in \K, 
and then $F$ preserves $\alpha$-filtered colimits.

\begin{remark}\label{rmk:accessible}
There is a corresponding result for the case where \K is accessible,
but we have not taken the trouble to formulate it here, since as 
usual there is a greater sensitivity to the choice of regular cardinal
in the accessible case than in the locally presentable one.
\end{remark}

The totality of small functors from $\K\op$ to \V forms a \V-category \PK
which is cocomplete and is in fact the free cocompletion of \K via the
Yoneda embedding $Y:\K\to\PK$. In the case where \K is small, \PK is
simply the presheaf category $[\K\op,\V]$, but in general not every
presheaf is small.

\begin{example}\label{large}
Let \V be \Set, and let \K be any large set $X$, seen as a discrete
category. Then a presheaf on \K can be seen as an $X$-indexed 
set $A\to X$, and it is small if and only if $A$ is so.
\end{example}

The construction \PK is pseudofunctorial in \K, and forms part of 
a pseudomonad \P on \VCat. We shall also consider free completions 
under certain types of colimit. Let $\Phi$ be a class of weights
with small domain. For a \V-category \K write $\Phi(\K)$ for
the closure of \K in \PK under $\Phi$-colimits. The Yoneda embedding
$Y:\K\to\Phi(\K)$ exhibits $\Phi(\K)$ as the free completion of
\K under $\Phi$-colimits. The class $\Phi$ is said to be {\em saturated}
if, whenever \K is small, $\Phi(\K)$ consists exactly of the presheaves
on \K lying in $\Phi$. (This idea goes back to \cite{Albert-Kelly}, where
the word ``closed'' was used rather than ``saturated''.)
Once again the construction $\Phi(\K)$  is
pseudofunctorial in \K and forms part of a pseudomonad $\Phi^*$ on \VCat.
The union $\Phi^*$ of all the $\Phi(\C)$ with \C small is a new 
class of weights called the {\em saturation} of $\Phi$. 

Thus far we have spoken only of smallness of presheaves, but we shall
also have cause to consider smallness of more general functors. Once
again, we say that a \V-functor $S:\K\to\M$ is small if it is the 
left Kan extension of some \V-functor $\C\to\M$ with small domain,
or equivalently, if it is the left Kan extension of its restriction
to some small full subcategory of \K. This definition works best when
\M is cocomplete, so that one can form the relevant left Kan extensions,
and we shall only use it in this context.  An important case is 
where $\M=[\C,\V]$ for some small \C. We say that $S:\K\to[\C,\V]$ is
{\em pointwise small} if the composite of $S$ with each evaluation functor
$\ev_C:[\C,\V]\to\V$ is small.

\begin{lemma}
A functor $S:\K\to[\C,\V]$ is small if and only if it is pointwise
small.
\end{lemma}

\proof
Since the evaluaton functors preserve Kan extensions the ``only if'' part
is immediate. Conversely, if $S$ is pointwise small, then for each
$C$ there is a small full subcategory $\B_C$ of \K with the property
that $\ev_C S$ is the left Kan extension of its restriction to $\B_C$.
Since \C is small, the union \B of the $\B_C$ is small, and now each
$\ev_C S$ is the left Kan extension of its restriction to \B, hence
the same is true of $S$.
\endproof

In Section~\ref{sect:M} we shall also consider the case where 
\M is locally presentable.
 
%%%%%%%%%%%%%%%%%%%%%%%%%%%%%%%%%%%%%%%%%%%%%%%%%%%%%%%%%%%%%%%%%%%%%%%%

\section{Limits of small functors}
\label{sect:limits}

%%%%%%%%%%%%%%%%%%%%%%%%%%%%%%%%%%%%%%%%%%%%%%%%%%%%%%%%%%%%%%%%%%%%%%%%

As observed above, if \K is small then \PK is the full presheaf category
$[\K\op,\Set]$ which is of course not just 
cocomplete but also complete. In general, however, a category of the form 
\PK need not be complete, as the following example, based on 
Example~\ref{large} shows:

\begin{example}\label{counterex}
If \V is \Set and \K is a large discrete category then \PK has
no terminal object.
%Let \V be \Set, and let \K be any large set $X$, seen as a discrete 
%category. Then a presheaf on \K is an $X$-indexed set $A\to X$, and it
%is small if and only if $A$ is so. Clearly \PK has no terminal object. 
\end{example}

We investigate which categories \K have the property that \PK is
complete. First observe that since \PK contains the representables,
any limit in \PK must be formed pointwise. Thus the question ``is \PK
complete?'' may be rephrased as ``are limits of small presheaves
on \K small?'' This may appear to involve consideration of the illegitimate
$[\K\op,\V]$, but in fact this is unnecessary. Given a weight 
$\phi:\C\to\V$, where \C is small, and a diagram 
$S:\C\to\PK$, we may regard $S$ as a functor $\bar{S}:\K\op\to[\C,\V]$,
and compose $\bar{S}$ with $\{\phi,-\}:[\C,\V]\to\V$, and ask whether the
composite $\{\phi,\bar{S}-\}:\K\op\to\V$ is small.

An arbitrary $\bar{S}:\K\op\to[\C,\V]$ arises in this way from some
$S:\C\to\PK$ if and only if $\bar{S}$ is {\em pointwise small}; recall
from the previous section that this means that each $\ev_C\bar{S}:\K\op\to\V$
is small, but that it is equivalent to $\bar{S}$ itself being small.

\begin{proposition}
The limit of $S:\C\to\PK$ weighted by $\phi:\C\to\V$ exists if
and only if $\{\phi,\bar{S}-\}$ is small; \PK has all $\phi$-limits
if and only if $\{\phi,R-\}$ is small for every small
$R:\K\op\to[\C,\V]$.
\end{proposition}

Related to the existence of limits in \K is the existence of 
a right adjoint to $\P F:\PK\to\PL$ for a functor $F:\K\to\LL$.
Here $\P F$ is given by left Kan extensions along $F$, so if
\K were small then $\P F$ would have a right adjoint given by
restriction along $F$. In general, however, the restriction $GF$
of a small $G:\LL\op\to\V$ need not be small; indeed the restriction
$\LL(F,L):\K\op\to\V$ of a representable $\LL(-,L)$ need not be small.
But if each $\LL(F,L)$ is small, we have the right adjoint:
%If \K is small, then restriction along $F$ gives such a right
%adjoint, but in general the restriction $GF$ of a small functor
%$G:\LL\op\to\V$ need not be small. However we have:

\begin{proposition}\label{prop:restriction}
For an arbitrary functor $F:\K\to\LL$, there is a right adjoint
to $\P F:\PK\to\PL$ if and only if $\LL(F,L):\K\op\to\V$ is small
for every object $L$ of \LL, and then the right adjont is given
by restriction along $F$.
\end{proposition}

\proof
If $\P F$ has a right adjoint $R$, then 
$$RGA\cong\PK(YA,RG)\cong\PL(\P F.YA,G)\cong\PL(YFA,G)\cong GFA$$
for any $G$ in \PL, and so $R$ must be given by restriction along
$F$. Thus $RYL=\LL(F,L)$, which must therefore be small. 

Suppose conversely that each $\LL(F,L)$ is small. Each $G$ in \PL
is a small colimit of representables. Since restricting along $F$
preserves colimits, $GF$ is a small colimit of functors of
the form $\LL(F,L)$, but these are small by assumption, so $GF$ is
small.
\endproof

Our first example of a large category \K with \PK complete is
the opposite of a locally presentable category.

\begin{proposition}
\PK is complete if $\K\op$ is locally presentable.
\end{proposition}

\proof
If $\K\op$ is locally presentable and $R:\K\op\to[\C,\V]$ is  small,
then for each object $C$ of \C there is a regular cardinal $\alpha_C$
for which $\ev_C R$ is $\alpha_C$-accessible. Since \C is small, we
may choose a regular cardinal $\alpha$ for which $\K\op$ is an 
$\alpha$-accessible category, $R$ is an $\alpha$-accessible functor,
and $\phi$ is $\alpha$-presentable in $[\C,\V]$. Then $R$ and $\{\phi,-\}$ 
preserve $\alpha$-filtered colimits, hence so does $\{\phi,R-\}$.
\endproof

\begin{remark}
The proposition remains true if $\K\op$ is accessible; the
comments made in Remark~\ref{rmk:accessible} still apply.
\end{remark}

\begin{corollary}\label{cor:smallA}
\PK is complete if \K is $[\A,\V]\op$ for a small category \A.
\end{corollary}

In other words, \PK is complete if $\K=\P(\A\op)\op$ for a small \A.
We shall now show how to remove the hypothesis that \A is small.
First observe $\P J:\PK\to\PL$ is given by left Kan extension along
$J$, so is fully faithful if $J$ is so.

\begin{proposition}\label{prop:PLimK}
\PK is complete if $\K=\P(\LL\op)\op$.
\end{proposition}

\proof Let \C be a small category and let $R:\K\op\to[\C,\V]$ be
small; we must show that $\{\phi,R-\}$ is small. Now $R$ is
the left Kan extension of its restriction to a small full subcategory
\D of $\P(\LL\op)$. Each $D\in\D$ is a small functor $\LL\to\V$, so
is the left Kan extension of its restriction to some small $\B_D$. 
The union \B of the $\B_D$ is small, and now the full inclusion 
$J:\B\op\to\LL\op$ induces a full inclusion $\P J:\P(\B\op)\to\P(\LL\op)$ 
containing \D.

Now \B is small, so $\P J$ has a right adjoint $J^*$ given by restriction
along $J$, and thus $\Lan_{\P J}$ is itself given by restriction along $J^*$.
Since $R$ is the left Kan extension of its restriction $S$ along $\P J$,
we have
$$\{\phi,R-\} = \{\phi,-\}R \cong \{\phi,-\} \Lan_{\P J} S \cong 
\{\phi,-\} SJ^* \cong \Lan_{\P J} \{\phi,-\}S = \Lan_{\P J}\{\phi,S-\}$$
and so $\{\phi,R-\}$ will be small if $\{\phi,S-\}$ is so.
Now $S:\P(\B\op)\to\V$ is the left Kan extension of its restriction to \D,
hence small, and \B is small, so by Corollary~\ref{cor:smallA} we conclude
that $\{\phi,S-\}$ is small.
\endproof

We are now ready to prove the main result of this section:

\begin{theorem}\label{thm:main}
\PK is complete if and only if it has limits of representables.
\end{theorem}

\proof
The ``only if'' part is trivial, so suppose that \PK has limits
of representables. Let $\LL=\P(\K\op)\op$, and let $Z:\K\to\LL$ be
the Yoneda embedding. By Proposition~\ref{prop:restriction} the
fully faithful $\P Z:\PK\to\PL$ has a right adjoint if $\LL(Z,L)$
is small for each $L$. But $\LL(Z,L)=\P(\K\op)(L,Y)$, where
$L:\K\to\V$ is a small functor. Then $L$ is the left Kan extension
of its restriction to some small full subcategory $J:\B\to\K$, 
and now $\P(\K\op)(\Lan_J(LJ),Y)=\P(\B\op)(LJ,YJ)$ which is the 
$LJ$-weighted limit of a diagram of representables, thus small
by assumption. This proves that \PK is a full coreflective subcategory
of \PL; since \PL is complete by Proposition~\ref{prop:PLimK}, it 
follows that \PK is so.
\endproof

\begin{corollary}\label{cor:complete}
\PK is complete if \K is so.
\end{corollary}

%%%%%%%%%%%%%%%%%%%%%%%%%%%%%%%%%%%%%%%%%%%%%%%%%%%%%%%%%%%%%%%%%%%%%%%%%%%

\section{Particular types of limit}
\label{sect:sound}

%%%%%%%%%%%%%%%%%%%%%%%%%%%%%%%%%%%%%%%%%%%%%%%%%%%%%%%%%%%%%%%%%%%%%%%%%%%

This section gives a more refined result, dealing with particular classes
of limits. It also provides an alternative proof for the main results of
the previous section. It is based on the ideas of \cite{ABLR}.

Let $\Phi$ be a class of weights. For a \V-category \C, we write
$\Phi\C$ for the closure of the representables in \PC under 
$\Phi$-weighted colimits. We suppose that the class $\Phi$ satisfies
the following conditions:
\begin{enumerate}[(a)]
\item (smallness) If \C is small then so is $\Phi\C$;
\item (soundness) If \D is small and $\Phi$-complete, and $\psi:\D\to\V$
is $\Phi$-continuous, then $\psi*-:[\D\op,\V]\to\V$ is $\Phi$-continuous.
\end{enumerate}
%Write $\Phi^*$ for the saturation of $\Phi$: a presheaf $\psi:\C\op\to\V$
%is in $\Phi^*$ if and only if it is in the closure in $[\C\op,\V]$ of
%the representables under $\Phi$-colimits.

\begin{example}~
\begin{enumerate}
\item If \V is \Set, then any sound doctrine in the sense of \cite{ABLR}
provides an example. Thus one could take $\Phi$ to be the (class of weights
corresponding to the) finite limits, or the $\alpha$-small limits
for some regular cardinal $\alpha$, or the finite products,
or the finite connected limits.
\item For any locally $\alpha$-presentable \V, by the results of 
\cite[(6.11),(7.4)]{Kelly-amiens} 
one can take 
$\Phi$ to be the class $\P_\alpha$ of $\alpha$-small limits.
\item If \V is cartesian closed, then by the results of \cite{Borceux-Day} 
(see also \cite{fpp}) one can take $\Phi$ to be the class of finite products. 
In fact by the results of \cite{Borceux-Day} this is still the case if \V is 
the algebras of any commutative finitary theory over a cartesian closed
category.
\end{enumerate}
\end{example}

\begin{lemma}
If \K is $\Phi$-cocomplete and $J:\C\to\K$ is a small full subcategory,
then the closure $\bar{\C}$ of \C in \K under $\Phi$-colimits is small.
\end{lemma} 

\proof
By the smallness assumption on $\Phi$, the free $\Phi$-cocompletion 
$\Phi\C$ of \C is small. Then $\bar{\C}$ is given, up to equivalence,
by the full image of the $\Phi$-cocontinuous extension $\bar{J}:\Phi\C\to\K$
of $J$; thus $\bar{\C}$ is small since $\Phi\C$ is so.
\endproof

\begin{proposition}
If \K is $\Phi$-complete then so is $\PK$. 
\end{proposition}

\proof
Let $\phi:\C\to\V$ be in $\Phi$, with \C small, and let $S:\K\op\to[\C,\V]$
be small. Then $S$ is the left Kan extension of its restriction to some small
full subcategory $J\op:\B\op\to\K\op$. By the lemma, we may suppose \B to be
closed in \K under $\Phi$-limits. Thus $S=\Lan_{J\op}R$, where \B is
small and $\Phi$-complete, $J:\B\to\K$ is $\Phi$-continuous, and 
$R:\B\op\to[\C,\V]$. Now $\K(K,J-):\B\to\V$ is $\Phi$-continuous for all
$K\in\K$, so $\K(K,J)*-:[\B\op,\V]\to\V$ is $\Phi$-continuous, so 
$\{\phi,-\}:[\C,\V]\to\V$ preserves the left Kan extension $S=\Lan_{J\op}R$.
In other words
$$\{\phi,S\}=\{\phi,\Lan_{J\op}R\}=\Lan_{J\op}\{\phi,R\}$$
and so $\{\phi,S\}$ is small.
\endproof

\begin{proposition}
\PK has all $\Phi$-limits if and only if it has $\Phi^*$-limits of 
representables.
\end{proposition}

\proof
Recall from Section~\ref{sect:small} that $\Phi^*$ is the ``saturation''
of $\Phi$, so that a \V-category has $\Phi$-limits if and only if
it has $\Phi^*$-limits, and in particular if \PK has $\Phi$-limits
then it certainly has $\Phi^*$-limits of representables. 

Let $Z:\K\to\LL$ be the free $\Phi$-completion of \K under $\Phi$-limits;
explicitly, $\LL=\Phi^*(\K\op)\op$, and $Z$ is the restricted Yoneda
embedding.
Then $\PL$ is $\Phi$-complete by the previous proposition. Since 
$\P Z:\PK\to\PL$ is fully faithful, $\PK$ will be $\Phi$-complete 
provided that $\P Z$ has a right adjoint. But this will happen if and
only if $\LL(Z-,F):\K\op\to\V$ is small for all $F\in\LL$. Now
$$\LL(Z-,F)=\P(\K\op)(F,Y-)$$
and the latter is an $F$-weighted limit of representables, with $F\in\Phi^*$.
\endproof

This provides an alternative proof of:

\begin{corollary}
\PK is complete if \K is.
\end{corollary}

\proof
If \K is complete, then it is $\P_{\alpha}$-complete for any regular
cardinal $\alpha$. Thus $\PK$ is $\P_{\alpha}$-complete for any regular
cardinal $\alpha$, and so is complete.
\endproof

%Notice that some restriction on $\Phi$ is necessary:

%\begin{example}
%Let $\V=\Set$, and let $\Phi$ consist of (the weight for) pullbacks.

%%%%%%%%%%%%%%%%%%%%%%%%%%%%%%%%%%%%%%%%%%%%%%%%%%%%%%%%%%%%%%%%%%%%%%%%%%%

\section{The case where $\V_0$ is a presheaf category}
\label{sect:ordinary}

%%%%%%%%%%%%%%%%%%%%%%%%%%%%%%%%%%%%%%%%%%%%%%%%%%%%%%%%%%%%%%%%%%%%%%%%%%%

For the first part of this section we suppose that $\V=\Set$, leading
to Theorem~\ref{thm:Set}. The latter should be attributed to Freyd,
although it may not have been written down by him in exactly this form;
it is a special case of \cite[Theorem~4.8]{Karazeris-completeness}. We include 
it as a warm-up for the more general case where the underlying
category $\V_0$ of \V is a presheaf category. This includes the case of 
the cartesian closed categories of directed graphs, or of simplicial sets, as
well as such non-cartesian cases as the category of $G$-graded sets,
for a group $G$, or the category of $M$-sets, for a commutative monoid $M$.

Suppose then that $\V=\Set$. First
observe that the statement {\em \PK has limits if and only if it has limits 
of representables} remains true if by limit we mean conical limit. To say
that \PK has conical limits of representables is to say that for
any $S:\C\to\K$ with \C small, the limit of $YS$ is small. But the
limit of $YS$ is the functor $\cone(S):\K\op\to\Set$ sending an object
$A$ to the set of all cones under $S$ with vertex $A$. To say that
this functor is small, is to say that there is a small full subcategory
\B of \K for which (i) any cone $\alpha:\Delta A\to S$ factorizes through
a cone $\beta:\Delta B\to S$ with $B\in\B$, and (ii) given cones
$\beta:\Delta B\to S$ and $\beta':\Delta B'\to S$ with $B,B'\in\B$, and 
arrows $f:A\to B$ and $f':A\to B'$ with $\beta.\Delta f=\beta'.\Delta f'$,
there is a ``zigzag'' of cones from $\beta$ to $\beta'$ with vertices in
\B. 

The existence of a small full subcategory \B satisfying (i) is clearly 
equivalent to the existence of a small set of cones through which every
cone factorizes: this is the {\em solution set condition}. In fact, however, 
if this solution set condition holds for any $S:\D\to\K$ with \D small
then \PK is complete, for we shall show below that if
\B satisfies (i), then we may enlarge \B to a new small full subcategory
\BB which satisfies (i) and (ii).  This is done as follows. Let $\B_0=\B$.
We construct inductively small full subcategories $\B_n$ for each natural
number $n$, and then define \BB to be the union of the $\B_n$. 

Let \DD be the category obtained from \D by freely adjoining two cones, 
with vertices 0 and 1, say. Let $\B_n$ be a small full subcategory,
and consider all functors $S':\DD\to\K$ extending \D, and sending 
the vertices 0 and 1 to objects of $\B_n$. For each such $S'$, we may 
by hypothesis choose a small full subcategory $\B_{S'}$ of \K which is 
a ``solution set'' for $S'$. Take $\B_{n+1}$ to be the union of $\B_n$
and all the $\B_{S'}$. This is a small union of small full subcategories,
so is itself a small full subcategory. Once again, \BB is a small union of the 
small full subcategories $\B_n$, and so is small. Clearly it satisfies (i); we 
check that it satisfies (ii) as well. Suppose then that $\beta:\Delta B\to S$
and $\beta':\Delta B'\to S$ are cones over $S$ with $B,B'\in\BB$, and that
$f:A\to B$ and $f':A\to B'$ are arrows with $\beta.\Delta f=\beta'.\Delta f'$.
Then $B$, $B'$, $\beta$, and $\beta'$ together define a functor $S':\DD\to\K$
extending $S$; while to give $A$, $f$, and $f'$ is precisely to give a cone
under $S'$. Since $B,B'\in\BB$, there is some $n\in\N$ for which 
$B,B'\in\B_n$, so there is a cone $S'$ with vertex $C$ in $\B_{S'}$ through 
which the cone $(A,f,f')$ factorizes. But this $C$ is in $\B_{n+1}$, and 
so in \BB. This proves:

\begin{theorem}\label{thm:Set}
\PK is complete if and only if for every diagram $S:\D\to\K$ with \D
small, there is a small set of cones under $S$ through which every 
cone factorizes.
\end{theorem}

We now extend this argument to the case where $\V_0$ is a presheaf
category $[\G\op,\Set]$. 
To extend the argument, we need to assume that the \V-category
\K admits tensors and cotensors by the representables; this means that for 
all $A,B\in\K$ and $G\in\G$, there are natural isomorphisms
$$\K(A,G\ct B) \cong [\G\op,\Set](\G(-,G),\K(A,B)) \cong \K(G\cdot A,B) $$
for objects $G\cdot A$ and $G\ct B$ of \K; the first operation is 
called a  {\em tensor} by $G$ and the second a {\em cotensor} by $G$.
When these exist, we say that \K is \G-tensored and \G-cotensored.
(Of course in the case $\V=\Set$ we have $\G=\{I\}$ and so this is automatic.) 

\begin{proposition}
If \K is \G-cotensored, then \PK is complete if and only if its underlying
ordinary category $(\PK)_0$ has conical limits of representables.
\end{proposition}

\proof
Recall \cite[3.10]{Kelly-book} that every weighted limit has a canonical expression 
as a conical limit of cotensors. On the other hand, since every object
of $\V_0=[\G\op,\Set]$ is (canonically) a conical colimit of representables,
and we have $(\colim_i G_i)\ct A\cong\lim_i(G_i\ct A)$, it follows that
every cotensor is canonically a conical limit of \G-cotensors,
and so finally that every weighted limit is canonically a conical limit
of \G-cotensors. Suppose now that we have a diagram of
representables $YS:\D\to\PK$ and a weight $\phi$. We can therefore express
this as a conical limit of \G-cotensors of representables. But \K was
assumed to be \G-cotensored, so a \G-cotensor of representables in \PK
exists, and is representable. Thus \PK will have weighted limits of 
representables provided that it has conical limits of representables. 
Finally, \PK is cocomplete, so is certainly tensored, thus
conical limits in \PK exist provided that they exist in the underlying
ordinary category $(\PK)_0$ of \PK, consisting of small \V-functors
$\K\op\to\V$ and \V-natural transformations between them.
\endproof

We now adapt the argument from the $\V=\Set$ case to the 
$V=[\G\op,\Set]$ case to prove:

\begin{proposition}
If \K is \G-tensored, then $(\PK)_0$ has conical limits of representables
if and only if, for every diagram $S:\C\to\K_0$, there is a small set
of cones through which every cone factorizes.
\end{proposition}

\proof
Suppose that $(\PK)_0$ has conical limits of representables, and let
$S:\C\to\K_0$ be given. Then $YS:\C\to(\PK)_0$ has a limit $L$ with
cone $\eta_C:L\to\K(-,SC)$. Also $L$ is a small colimit of representables
$\colim_i\K(-,B_i)$. For each $i$, there is an induced cone 
$\K(-,B_i)\to\K(-,SC)$, or equivalently $\beta_i C:B_i\to SC$ under $S$.
We claim that any cone $\alpha C:A\to SC$ factorizes through one of these.
Now $\K(-,\alpha C):\K(-,A)\to\K(-,SC)$ must factorize through $L$, but 
$\K(-,A)$ is representable, so homming out of it preserves colimits, and
so we get a factorization $\K(-,A)\to\K(-,B_i)$ for some $i$, and so
the desired $A\to B_i$.

For the harder part, suppose that for each $S:\C\to\K_0$, there is a
small set of cones $\beta_i C:B_i\to SC$ through which each cone 
factorizes. Then there is a small full
subcategory $\B_S$ of \K such that each cone under $S$ factorizes through
one whose vertex is in $\B_S$. Exactly as before, we construct
the (possibly larger but still) small full subcategory $J:\BB\to\K$ with
the property that if two cones with vertices in \BB are connected, then 
they are connected using cones with vertices in \BB. This implies that
for all $A\in\K$, we have $\lim_C \K_0(A,SC)\cong \colim_{B\in\BB}\K_0(A,B)$.
For any $G\in\G$, we have
\begin{align*}
\V_0(G,\lim_C\K(A,SC)) &\cong \lim_C\V_0(G,\K(A,SC)) \\
                       &\cong \lim_C\K_0(G\cdot A,SC) \\
                       &\cong \colim_{B\in\BB}\K_0(G\cdot A,B) \\
                       &\cong \colim_{B\in\BB}\V_0(G,\K(A,B)) \\
                       &\cong \V_0(G,\colim_{B\in\BB}\K(A,B))
\end{align*}
where the last step uses the fact that $G$ is representable, so $\V_0(G,-)$
preserves colimits. Now \G is dense in $\V_0$, so we have
\begin{align*}
\lim_C\K(A,SC) & \cong \colim_{B\in\BB}\K(A,B) \\
\lim_C\K(-,SC) & \cong \colim_{B\in\BB}\K(-,B)
\end{align*}
but the left hand side is the presheaf $\K\op\to\V$ which is the pointwise
limit of $YS$, and which we are to prove small, while the right hand side
is a small colimit of representables, since \BB is small.
\endproof

Combining the last two results, we have:

\begin{theorem}
Suppose the underlying category $\V_0$ of \V is a presheaf category
$[\G\op,\Set]$ and that \K is a \G-tensored and \G-cotensored \V-category. 
Then \PK is complete if and only if the following condition is satisfied.
For every small ordinary category \C and every functor $S:\C\to\K$,
there is a small set of cones $\lambda C:B\to SC$ through which every 
such cone factorizes.
\end{theorem}

%%%%%%%%%%%%%%%%%%%%%%%%%%%%%%%%%%%%%%%%%%%%%%%%%%%%%%%%%%%%%%%%%%%%%%%%%%%

\section{Preservation of limits}
\label{sect:continuity}

%%%%%%%%%%%%%%%%%%%%%%%%%%%%%%%%%%%%%%%%%%%%%%%%%%%%%%%%%%%%%%%%%%%%%%%%%%%

Having studied the categories
\K for which \PK is complete, we now turn to the functors $F:\K\to\LL$
for which $\P F$ is continuous. 

We saw \PK is {\em always} complete if \K is small; the situation for functors
is totally different:

\begin{example}
Let \V=\Set, let \K be the terminal category 1, let \LL be the 
discrete category 2, and let $F:\K\to\LL$ be the first injection.
Then $\PK=\Set$ and $\PL=\Set^2$, which are of course complete;
but $\P F:\Set\to\Set^2$ is the functor sending a set $X$ to $(X,0)$,
which clearly fails to preserve the terminal object.
\end{example}

Consider a functor $F:\K\to\LL$, where \PK and \PL are complete, a
weight $\phi:\C\to\V$ and a diagram $S:\C\to\PK$. Let $R:\K\op\to[\C,\V]$
be the corresponding pointwise small functor. To say that $\P F$
preserves the limit $\{\phi,S\}$ is to say that $\{\phi,-\}:[\C,\V]\to\V$
preserves the left Kan extension $\Lan_F R$; that is, the colimit
$\LL(L,F-)*R$ for each object $L$ of \LL. 

\begin{proposition}
If $F:\K\to\LL$ is a right adjoint then $\P F$ is continuous.
\end{proposition}

\proof
If $F$ has a left adjoint $G$, then 
$$\{\phi,\LL(L,F-)*R\} \cong \{\phi,\K(GL,-)*R\} \cong \{\phi,RGL\}$$
while
$$\LL(L,F-)*\{\phi,R\} \cong \LL(GL,-)*\{\phi,R\} \cong \{\phi,RGL\}.$$
\endproof

Along the same lines, observe that $\P F.Y\cong YF$, so that if $\P F$
is continuous then $F$ must preserve any limits which exist.

Suppose that $F:\K\to\LL$ is given, with \PK and \PL complete.
Then $\P F$ is continuous if and only if each $\ev_L.\P F$ is so;
but $\ev_L.\P F$ is just $\LL(L,F)*-$. If \K is small, then 
$\LL(L,F)*-$ is continuous if and only if it is $\alpha$-continuous
for every regular cardinal $\alpha$; in other words, if 
$\LL(L,F)$ is $\alpha$-flat for every $\alpha$.

More generally, if \PK is $\alpha$-complete, we say that a functor $G:\K\to\V$
is $\alpha$-flat if $G*-:\PK\to\V$ is $\alpha$-continuous, and 
$\infty$-flat if $G*-$ is continuous; that is, if $G$ is $\alpha$-flat
for every $\alpha$. Thus $\P F$ will be continuous if and only if
each $\LL(L,F)$ is $\infty$-flat. Similarly, if $\Phi$ is a class
of weights satisfying the conditions in Section~\ref{sect:sound}
and \PK is $\Phi$-complete, we say that $G:\K\to\V$ is $\Phi$-flat
if $G*-$ is $\Phi$-continuous.

\begin{lemma}
If \K is complete and $G:\K\to\V$ continuous then $\Lan_Y G:\PK\to\V$
is continuous.
\end{lemma}

\proof
First observe that if \K is complete then \PK is so. Let $\phi:\D\to\V$
be a weight, and let $S:\D\to\PK$ correspond to the pointwise small
functor $R:\K\op\to[\D,\V]$. For $X\in\PK$, we have the formula
$(\Lan_Y G)X=X*G$, thus to say that $\Lan_Y G$ preserves the
limit $\{\phi,S\}$ is to say that $\{\phi,-\}$ preserves the colimit
$G*R$. Since $R$ is pointwise small,
it is the left Kan extension of its restriction to some full subcategory
$\B\op$ of $\K\op$. Let $\alpha$ be a regular cardinal for which $\phi$
is $\alpha$-small. We may choose \B to be closed in \K under
$\alpha$-limits, then the inclusion $J:\B\to\K$ preserves 
$\alpha$-limits. Then
$G*R\cong G*(\Lan_J(RJ))\cong GJ*RJ$ by \cite[4.1]{Kelly-book}, 
and $GJ:\B\to\V$ preserves $\alpha$-limits, hence so does $GJ*-$ 
by \cite[(6.11),(7.4)]{Kelly-amiens}, and now
$$\{\phi,G*R\} \cong \{\phi,GJ*RJ\} \cong GJ*\{\phi,RJ\}\cong
G*\Lan_J\{\phi,RJ\}.$$
On the other hand $\Lan_J$ preserves $\alpha$-limits since $J$ does
so, thus
$$G*\Lan_J\{\phi,RJ\} \cong G*\{\phi,\Lan_J(RJ)\} \cong
G*\{\phi,R\}.$$
This proves that $G*-$ preserves the limit $\{\phi,R\}$, and
so that $\Lan_Y G$ preserves $\{\phi,S\}$.
\endproof

\begin{theorem}
Let \K and \LL be complete. Then $F:\K\to\LL$ is continuous if and 
only if $\P F:\PK\to\PL$ is so.
\end{theorem}

\proof
The ``if part'' was observed above. Suppose then that $F$ is continuous.
Then each $\LL(L,F)$ is continuous, so $\Lan_Y \LL(L,F)$ is continuous,
but $\Lan_Y \LL(L,F)\cong\ev_L.\P F$, and so $\P F$ is continuous, since
limits in \PL are constructed pointwise.
\endproof

\begin{remark}
The Yoneda embedding $Y:\K\to\PK$ preserves any existing limits,
and is continuous if \K is complete. The pseudomonad \P is of the
Kock-Z\"oberlein type, and so the multiplication $\P\!\PK\to\PK$
has both adjoints so also preserves any existing limits (or colimits).
Thus the pseudomonad \P lifts from \VCat to the 2-category of complete
\V-categories, continuous \V-functors, and \V-natural transformations.
\end{remark}

\begin{remark}
Suppose once again that $\Phi$ is a class of weights satisfying the
conditions of Section~\ref{sect:sound}.
Suppose that \K and \LL are $\Phi$-complete and $F:\K\to\LL$ is
$\Phi$-continuous. Then each $\LL(L,F)$ is $\Phi$-continuous, 
so each $\LL(L,F)*-$ is $\Phi$-continuous; that is, each
$\ev_L\P(F)$ is $\Phi$-continuous, and so finally $\P F:\P K\to\P L$
is $\Phi$-continuous. Thus the pseudomonad $\Phi^*$ lifts
from \VCat to the 2-category of $\Phi$-complete \V-categories, 
$\Phi$-continuous \V-functors, and \V-natural transformations.
\end{remark}

%%%%%%%%%%%%%%%%%%%%%%%%%%%%%%%%%%%%%%%%%%%%%%%%%%%%%%%%%%%%%%%%%%%%%%%%%%%

\section{Monoidal structure on \PK}
\label{sect:promonoidal}

%%%%%%%%%%%%%%%%%%%%%%%%%%%%%%%%%%%%%%%%%%%%%%%%%%%%%%%%%%%%%%%%%%%%%%%%%%%

In this section we suppose that \K is a \V-category for 
which \PK is complete. If \K is small, so that \PK is $[\K\op,\V]$,
monoidal closed structures on \PK correspond to promonoidal structures
on $\K\op$ \cite{Day-convolution}. These consist of \V-functors
$P:\K\op\ot\K\ot\K\to\V$ and $J:\K\op\to\V$ equipped with coherent
associativity and unit isomorphisms. 

If \K is large, we shall insist that $P(-;A,B):\K\op\to\V$ and $J:\K\op\to\V$
be small, and we write $P:\K\ot\K\to\PK:(A,B)\mapsto P(-;A,B)$ and $J\in\PK$. 
If $F,G\in\PK$ are
given, we define $F\ot G$ using the usual convolution formula:
$$F\ot G = \int^{A,B} P(-;A,B)\ot FA\ot GB.$$
This is small, since each $P(-;A,B)$ is small by assumption, so 
$\int^A P(-;A,B)\ot FA$ is a small ($F$-weighted) colimit of small 
presheaves for each $B$,
and so $\int^{A,B} P(-;A,B)\ot FA\ot GB$ is itself a small colimit of 
small presheaves, hence small.

In the usual case, where \K is small, this monoidal structure is 
closed, with (right) internal hom given by 
\begin{align*}
[G,H] &\cong \int_{B,C} [P(C;-,B)\ot GB,HC] \\
      &\cong \int_{B,C}[GB,[P(C;-,B),HC]]
\end{align*}
If \K is large, this need not lie in \PK, but if it does so, then 
it will still provide the internal hom. Now $G$ is small, and the
expression above for $[G,H]$ is precisely the $G$-weighted limit
of the functor sending $B$ to $\int_C[P(C;-,B),HC]$. Since \PK 
is complete this limit will exist provided that this functor actually
lands in \PK; that is, provided that 
$$\int_C[P(C;-,B),HC]:\K\op\to\V$$
is small for all $B\in\K$.

The case of the other internal hom is similar, and we have:

\begin{proposition}
The convolution monoidal category \PK is closed if and only if the
presheaves $\int_C[P(C;-,B),HC]$ and $\int_C[P(C;B,-),HC]$ are small
for all $B\in\K$.
\end{proposition}

An important special case is where the promonoidal structure $P$ is
a filtered colimit $P=\colim_i P_i$ of promonoidal structures $P_i$
which are in fact monoidal, as in 
$$P_i(C;A,B)=\K(C,A\ot_i B).$$
We call such a promonoidal structure $P$ {\em approximately monoidal}; 
of course every monoidal structure is approximately monoidal. (We are using
the fact that the colimit is filtered to obtain the associativity and unit
isomorphisms; a general colimit of promonoidal structures need not be
promonoidal.)

In the approximately monoidal case a simplification is possible, since
\begin{align*}
\int_C[P(C;-,B),HC] &\cong \int_C[\colim_i P_i(C;-,B),HC] \\
                    &\cong \lim_i\int_C[P_i(C;-,B),HC] \\
                    &\cong \lim_i\int_C[\K(C,-\ot_i B),HC] \\
                    &\cong \lim_i H(-\ot_i B)
\end{align*}
which is small provided each $H(-\ot_i B)$ is so. But $H$ is small,
so has the form $\Lan_J(HJ)$ for some $J:\D\to\K\op$ with \D small.
Then 
\begin{align*}
H(-\ot_i B) &\cong \int^D \K\op(D,-\ot_i B)\cdot HJD \\
            &\cong \int^D \K(-\ot_i B,D)\cdot HJD 
\end{align*}
which is a small ($HJ$-weighted) colimit of presheaves $\K(-\ot_i B,D)$ 
with $D\in\D$, so will be small provided that the $\K(-\ot_i B,D)$ are so. 
Once again the case of the other internal hom is similar, and we have:

\begin{proposition}
The convolution monoidal category \PK arising from an approximately 
monoidal structure on \K is closed if and only if the
presheaves $\K(-\ot_i B,D)$ and $\K(B\ot_i-,D)$ are small for all 
$B$ and $D$ in \K, and for each monoidal structure $\ot_i$.
\end{proposition}

In particular we have:

\begin{proposition}
The convolution monoidal category \PK arising from a
monoidal structure on \K is closed if and only if the
presheaves $\K(-\ot B,D)$ and $\K(B\ot-,D)$ are small for all 
$B$ and $D$ in \K.
\end{proposition}

\begin{example}\label{ex:Rosicky}~
\begin{enumerate}
\item The special case where $\V=\Set$ and the monoidal structure
is cartesian was proved in \cite{Rosicky-ccc}.
\item If \K is not just monoidal but closed then the $\K(-\ot B,D)$
and $\K(B\ot-,D)$ are not just small but representable, and so \PK 
is monoidal closed.
\item If \V is cartesian monoidal (so that $\ot=\t$), 
and $\K=\E\op$ where \E is also cartesian monoidal, then 
$\K(-\ot B,D)=\E(D,-\t B)=\E(D,B)\t\E(D,-)$ which is given by 
tensoring the representable $\E(D,-)$ by the \V-object $\E(D,B)$,
and so is small. Thus once again \PK is monoidal closed.
\end{enumerate}
\end{example}

%%%%%%%%%%%%%%%%%%%%%%%%%%%%%%%%%%%%%%%%%%%%%%%%%%%%%%%%%%%%%%%%%%%%%%%%%%%

\section{Functors with codomain other than \V}
\label{sect:M}

%%%%%%%%%%%%%%%%%%%%%%%%%%%%%%%%%%%%%%%%%%%%%%%%%%%%%%%%%%%%%%%%%%%%%%%%%%%

In this section we consider small functors $\K\op\to\M$ where \M is
cocomplete, building on our earlier work on the case $\M=\V$ and
$\M=[\C,\V]$. 

In that earlier work, we considered when, for a small functor
$S:\K\op\to[\C,\V]$, each $\{\phi,S\}$ was small. But $\{\phi,-\}$
is just the representable functor $[\C,\V](\phi,-)$, which motivates
the following definition: a functor $S:\K\to\M$ is {\em representably
small} if each $\M(M,S):\K\to\V$ is small. Thus Corollary~\ref{cor:complete}
asserts that if \K is complete then every small functor 
$\K\op\to[\C,\V]$ is representably small.

In this section we investigate the relationship between smallness and
representable smallness for more general \M. We have already seen 
that smallness does not in general imply representable smallness. 
For an explicit counterexample in the case $\M=\V$ we have:

\begin{example}
As in Example~\ref{counterex}
let \V be \Set, and let \K be any large set $X$, seen as a discrete 
category. Then a presheaf on \K is an $X$-indexed set $A\to X$, and it
is small if and only if $A$ is so. Certainly $x:1\to X$ is small, for 
any $x\in X$; this corresponds to the representable presheaf
$X(-,x):X\to\Set$ sending $x$ to 1, and all other elements to 0. Now
$\Set(0,X(-,x))$ is the terminal presheaf, which as we have seen is
not small. Thus $X(-,x)$ is small but not representably small.
\end{example}

To see that a representably small functor need not be small, we have:

\begin{example}
If \K is a large \V-category for which \PK is complete (for example if
\K is complete), then the Yoneda embedding $Y_{\K\op}:\K\op\to\P(\K\op)$
is representably small. For if $F\in\P(\K\op)$, the composite
$\P(\K\op)(F,Y):\K\op\to\V$ is 
the $F$-weighted limit of $Y:\K\to\PK$, so is small since $F$ is small 
and \PK is complete. But $Y_{\K\op}:\K\op\to\P(\K\op)$ is not small 
unless \K is so. For if $Y_{\K\op}$ were small, $\K\op$ would have a 
small full subcategory $J:\C\op\to\K\op$ for which $Y=\Lan_J(YJ)$, so
$$\K(-,A) = \int^{C\in\C} \K(JC,A)\cdot\K(-,JC)$$
for all $A$, and in particular 
$$\K(A,A)= \int^C \K(JC,A)\cdot\K(A,JC)$$
and so the identity $1:A\to A$ must factorize through some $JC$; in
other words, each $A\in\K$ is a retract of some object in \C. But this
clearly implies that \K is small.
\end{example}

As a first positive result we have:

\begin{proposition}
If \K is a \V-category for which \PK admits cotensors, a presheaf 
$F:\K\op\to\V$ is small if and only if it is representably small.
In particular this will be the case if \PK is complete. 
\end{proposition}

\proof
Representably small presheaves are always small, since $\V(I,F)$ is 
just $F$, for any presheaf $F$.
It remains to show that any small presheaf $F:\K\op\to\V$ is representably
small. Suppose then that $X\in\V$. Then $\V(X,F)$ is the cotensor 
$X\pitchfork F$ of $F$ by $X$, which is small by assumption.
\endproof

For the remainder of the section we suppose that \K is a \V-category
for which \PK is complete, and that \M is a locally presentable 
\V-category. If $\beta$ is a regular cardinal for which \M is locally 
$\beta$-presentable, write  $\M_\beta$ for the full subcategory of \M 
consisting of the $\beta$-presentable objects, and $W:\M\to[\M\op_\beta,\V]$ 
for the canonical (fully faithful) inclusion.

\begin{lemma}
For a \V-functor $S:\K\op\to\M$, the following are equivalent:
\begin{enumerate}[(a)]
\item $S$ is representably small;
\item $WS$ is small;
\item $S$ is small.
\end{enumerate}
\end{lemma}

\proof
$(a)\Rightarrow(b)$. To say that $S$ is representably small is to say
that $\M(M,S)$ is small for all $M\in\M$; to say that $WS$ is small is
to say that this is so for all $M\in\M_\beta$, so this is immediate.

$(b)\Rightarrow(c)$. For each $M\in\M_\beta$ we have $\M(M,S)$ small,
so it is the left Kan extension of its restriction to some full subcategory
$\D_M$ of $\K\op$. Since $\M_\beta$ is small, the union \D of the 
$\D_M$ is small, and each $\M(M,S)$ is the left Kan extension of its
restriction to \D. Thus $WS$ is the left Kan extension of its restriction
to \D. But $W$ is fully faithful, and so reflects Kan extensions; thus
also $S$ is the left Kan extension of its restriction to \D.

$(c)\Rightarrow(a)$. This is by far the hardest implication; we prove it
in several steps, analogous to the main steps used in 
preparation for the proof of Theorem~\ref{thm:main}.
Suppose then that $S$ is small and $M\in\M$; we must show that $\M(M,S)$ is
small. 

{\bf Case 1: $\K\op$ is locally presentable.} 
%For Case 1 only, replace
%(if necessary) $\beta$ by some greater regular cardinal for which 
%$\K\op$ is locally $\beta$-presentable.
Since $S$ is small, it is the left Kan extension $\Lan_{J\op}R$ along 
$J\op:\C\op\to\K\op$ of some $R:\C\op\to\M$ with \C small. Since \C
is small and $\K\op$ and \M are locally presentable, there exists a regular
cardinal $\gamma\ge\beta$ for which each $JC$ is $\gamma$-presentable
in $\K\op$ and $M$ is $\gamma$-presentable in \M. Now (a) $\K\op$ is
the free completion under $\gamma$-filtered colimits of the full
subcategory $(\K\op)_\gamma$ of $\K\op$ consisting of the $\gamma$-presentable
objects, (b) $S$ preserves $\gamma$-filtered colimits, and (c) $\M(M,-)$
preserves $\gamma$-filtered colimits. Thus $\M(M,S)$ preserves 
$\gamma$-filtered colimits, so is the left Kan extension of its restriction
to $(\K\op)_\gamma$. This proves that $\M(M,S)$ is
small, and so that $S$ is representably small.

{\bf Case 2: $\K\op=\P(\LL\op)$.} 
Then $S$ is the 
left Kan extension of its restriction to some small full subcategory
\D of $\P(\LL\op)$. Each $D\in\D$ is a small functor $\LL\to\V$, so is
the left Kan extension of its restriction to some small $\B_D$. The union
\B of the $\B_D$ is small, and now the full inclusion $J:\B\op\to\LL\op$
induces a full inclusion $\P J:\P(\B\op)\to\P(\LL\op)$ whose image
contains \D. 

Now \B is small, so $\P J$ has a right adjoint $J^*$ given by restriction
along $J$, and thus $\Lan_{\P J}$ is itself given by restriction along 
$J^*$. Since $S$ is the left Kan extension of its restriction $Q$ along
$\P J$, we have
$$\M(M,S) = \M(M,\Lan_{\P J}Q) = \M(M,QJ^*) = \M(M,Q)J^* = \Lan_{\P J}\M(M,Q)$$
and so $\M(M,S)$ will be small if $\M(M,Q)$ is so. Now $Q$ is the 
left Kan extension of its restriction to \D, hence small, so $\M(M,Q)$ is
small by Case 1. This proves that $\M(M,S)$ is small, and so that $S$ is
representably small.

{\bf Case 3: \PK is complete.}
The left Kan extension $\Lan_Y(S):\P(\K\op)\to\M$ of $S$ along the Yoneda 
embedding is small, so by Case~2 is representably small. Thus each 
$\M(M,\Lan_Y(S)):\P(\K\op)\op\to\V$ is
small; that is, a small colimit of representables. Now restriction 
along the Yoneda embedding preserves colimits, so it will send small
presheaves to small presheaves provided that it sends representables
to small presheaves; but the latter is equivalent to completeness of 
\PK. Thus each $\M(M,S)$ is small, and $S$ is representably small.
\endproof

Write $[\K\op,\M]_s$ for the \V-category of all small \V-functors
from $\K\op$ to \M.

\begin{theorem}
Let \M be a locally presentable \V-category, and $\K$ a \V-category 
for which \PK is complete. Then $[\K\op,\M]_s$ is complete.
\end{theorem}

\proof
%We shall show that limits in $[\K\op,\M]_s$ exist and are computed
%pointwise; to do this it suffices to show that the pointwise limit
%is in fact small, and so lies in $[\K\op,\M]_s$.
Let $\phi:\D\to\V$ and $S:\D\to[\K\op,\M]_s$ be given, where \D is small.
Since \D is small, the functor $\bar{S}:\K\op\to[\D,\M]$ corresponding to
$S$ is small. The ``pointwise limit'' is the composite
$$\xymatrix{\K\op \ar[r]^-{\bar{S}} & [\D,\M] \ar[r]^-{\{\phi,-\}} &
  \M}$$
and provided that this is small, and so lies in $[\K\op,\M]_s$, it will
be the limit. Since \M is locally presentable, by the lemma it will
suffice to show that each composite with $\M(M,-)$ is small. But for
any $X:\D\to\M$ we have 
\begin{align*}
\M(M,\{\phi,X\}) &\cong \{\phi,\M(M,X)\} \\
                 &\cong \int_D [\phi D,\M(M,XD)] \\
                 &\cong \int_D \M(\phi D\cdot M,XD) \\
                 &\cong [\D,\M](\phi_M,X)
\end{align*}
where $\phi_M:\D\to\M$ is the functor sending $D$ to $\phi D\cdot M$,
so now 
$\M(M,\{\phi,-\})$ is representable as 
$[\D,\M](\phi_M,-):[\D,\M]\to\V$.

Now $[\D,\M]$ is locally presentable, so by the lemma once again the small
$\bar{S}$ is representably small, and so $[\D,\M](\phi_M,\bar{S}):\K\op\to\V$
is small; but we have just seen that this is the composite of $\bar{S}$
with $\M(M,\{\phi,-\})$. This now proves that $\{\phi,-\}\circ\bar{S}$ is 
representably small, and so small, and it therefore provides the desired
limit $\{\phi,S\}$.
\endproof

%%%%%%%%%%%%%%%%%%%%%%%%%%%%%%%%%%%%%%%%%%%%%%%%%%%%%%%%%%%%%%%%%%%%%%%%%%%

\section{Isbell conjugacy}
\label{sect:Isbell}

%%%%%%%%%%%%%%%%%%%%%%%%%%%%%%%%%%%%%%%%%%%%%%%%%%%%%%%%%%%%%%%%%%%%%%%%%%%

If \C is a small category then as well as the Yoneda embedding 
$Y:\C\to[\C\op,\V]$ there is also the ``dual'' Yoneda embedding
$Z:\C\to[\C,\V]\op$, and this induces an adjunction between $[\C\op,\V]$
and $[\C,\V]\op$ called ``Isbell conjugacy''. The left adjoint 
$L:[\C\op,\V]\to[\C,\V]\op$ is given by $\Lan_Y Z$.

What happens if we replace \C be an arbitrary category \K? Then we 
have $Y:\K\to\PK$ and $Z:\K\to\P(\K\op)\op$, but do we still have
the adjunction between them? A sufficient condition for the left
adjoint $L:\PK\to\P(\K\op)\op$ to exist is that $\P(\K\op)\op$ be
cocomplete, or equivalently $\P(\K\op)$ complete, but in fact this
is also necessary. For if $\Lan_Y Z$ does exist, then for each 
small $F:\C\op\to\V$ the colimit $F*Z$ in $\P(\K\op)\op$ exists.
But then for any $\phi:\C\op\to\V$ and $S:\C\to\K$, we have
$\Lan_S \phi$ small, and $(\Lan_S \phi)*Z=\phi ZS$, and so $\P(\K\op)\op$
has arbitrary colimits of representables, $\P(\K\op)$ has arbitrary
limits of representables, and so $\P(\K\op)$ is in fact complete.

Thus $\Lan_Y Z:\PK\to\P(\K\op)\op$ exists if and only if $\P(\K\op)$
is complete, and dually the putative right adjoint 
$\P(\K\op)\op\to\PK$ exists if and only if \PK is complete. 

In particular, both will exist if \K is complete and cocomplete.

%%%%%%%%%%%%%%%%%%%%%%%%%%%%%%%%%%%%%%%%%%%%%%%%%%%%%%%%%%%%%%%%%%%%%%%%%%%

\bibliographystyle{plain}

%\bibliography{my}

\begin{thebibliography}{10}

\bibitem{ABLR}
Ji{\v{r}}{\'{\i}} Ad{\'a}mek, Francis Borceux, Stephen Lack, and
  Ji{\v{r}}\'\i\xspace Rosick{\'y}.
\newblock A classification of accessible categories.
\newblock {\em J. Pure Appl. Algebra}, 175(1-3):7--30, 2002.

\bibitem{AR}
Ji{\v{r}}{\'{\i}} Ad{\'a}mek and Ji{\v{r}}\'\i\xspace Rosick{\'y}.
\newblock {\em Locally presentable and accessible categories}, volume 189 of
  {\em London Mathematical Society Lecture Note Series}.
\newblock Cambridge University Press, Cambridge, 1994.

\bibitem{Albert-Kelly}
M.~H. Albert and G.~M. Kelly.
\newblock The closure of a class of colimits.
\newblock {\em J. Pure Appl. Algebra}, 51(1-2):1--17, 1988.

\bibitem{Borceux-Day}
Francis Borceux and B.~J. Day.
\newblock On product-preserving {K}an extensions.
\newblock {\em Bull. Austral. Math. Soc.}, 17(2):247--255, 1977.

\bibitem{Chorny}
B. Chorny and W.G. Dwyer.
\newblock {\em Homotopy theory of small diagrams over large categories}
\newblock Preprint.

\bibitem{Day-convolution}
Brian Day.
\newblock On closed categories of functors.
\newblock In {\em Reports of the Midwest Category Seminar, IV}, Lecture Notes
  in Mathematics, Vol. 137, pages 1--38. Springer, Berlin, 1970.

\bibitem{Freyd-lucid}
Peter Freyd.
\newblock Several new concepts: {L}ucid and concordant functors, pre-limits,
  pre-completeness, the continuous and concordant completions of categories.
\newblock In {\em Category Theory, Homology Theory and their Applications, III
  (Battelle Istitute Conference, Seattle, Wash., 1968, Vol. Three)}, pages
  196--241. Springer, Berlin, 1969.

\bibitem{Gabriel-Ulmer}
Peter Gabriel and Friedrich Ulmer.
\newblock {\em Lokal pr\"asentierbare {K}ategorien}.
\newblock Springer-Verlag, Berlin, 1971.

\bibitem{Karazeris-completeness}
Panagis Karazeris, Ji{\v{r}}{\'{\i}} Rosick{\'y}, and Ji{\v{r}}{\'\i} Velebil.
\newblock Completeness of cocompletions.
\newblock {\em J. Pure Appl. Algebra}, 196(2-3):229--250, 2005.

\bibitem{Kelly-amiens}
G.~M. Kelly.
\newblock Structures defined by finite limits in the enriched context. {I}.
\newblock {\em Cahiers Topologie G\'eom. Diff\'erentielle}, 23(1):3--42, 1982.

\bibitem{Kelly-book}
G.~M. Kelly.
\newblock Basic concepts of enriched category theory.
\newblock {\em Repr. Theory Appl. Categ.}, (10):vi+137 pp. (electronic), 2005.
\newblock Originally published as LMS Lecture Notes 64, 1982.

\bibitem{fpp}
G.~M. Kelly and Stephen Lack.
\newblock Finite-product-preserving functors, {K}an extensions and
  strongly-finitary {$2$}-monads.
\newblock {\em Appl. Categ. Structures}, 1(1):85--94, 1993.

\bibitem{vcat}
G.~M. Kelly and Stephen Lack.
\newblock $\mathscr{V}$-{C}at is locally presentable or locally bounded if
  $\mathscr{V}$ is so.
\newblock {\em Theory Appl. Categ.}, 8:555--575, 2001.

\bibitem{Lindner}
Harald Lindner.
\newblock Enriched categories and enriched modules.
\newblock {\em Cahiers Topologie G\'eom. Diff\'erentielle}, 22(2):161--174,
  1981.

\bibitem{Rosicky-ccc}
J.~Rosick{\'y}.
\newblock Cartesian closed exact completions.
\newblock {\em J. Pure Appl. Algebra}, 142(3):261--270, 1999.

\end{thebibliography}

\end{document}